\title{Computing Adem Cohomology Operations}
\author{
        Rocio Gonzalez--Diaz, Pedro Real\\
         \normalsize Depto. de Matematica Aplicada I, E.T.S.I. Informatica,\\ 
         \normalsize Universidad de
Sevilla (Spain)\\
\normalsize e-mail: $\{$rogodi,real$\}$@us.es
}
\date{July, 2006}

\documentclass[12pt]{article}

\begin{document}

\newcommand{\0}{\bar{0}}
\newcommand{\ra}{\rightarrow}
\newcommand{\la}{\leftarrow}
\newcommand{\Ra}{\Rightarrow}
\newcommand{\D}{\displaystyle}
\newcommand{\scst}{\scriptscriptstyle}
\newcommand{\scsz}{\scriptsize}
\newcommand{\pa}{\partial}
\newcommand{\m}{\bar{m}}
\newcommand{\ot}{\otimes}
\newcommand{\ti}{\times}
\newcommand{\al}{\alpha}
\newcommand{\be}{\beta}
\newcommand{\Z}{{\bf Z}}
\newcommand{\bb}{\bar{\be}}
\newcommand{\Zr}{{\bf Z}_2}
\newcommand{\tc}{{\theta\cap}}
\newcommand{\1}{\bar{1}}

\newtheorem{lem}{Lemma}
\newtheorem{thm}[lem]{Theorem}
\newtheorem{cor}[lem]{Corollary}
\newtheorem{prop}[lem]{Proposition}
\newtheorem{alg}[lem]{Algorithm}
\maketitle

\begin{abstract}
We deal with the problem of obtaining explicit simplicial formulae
defining the classical Adem cohomology operations at the cochain
level.
Having these formulae at hand, we design an algorithm for computing
these operations for any finite simplicial set.
\newline
$\mbox{ }$
\end{abstract}

\section{Introduction}

Cohomology operations are algebraic operations on the cohomology groups of spaces. This machinery
is useful when the (co)homology and the cup product on cohomology fail to distinguish two spaces.
In the literature, there exists classical methods for computing the (co)homology
of spaces (see \cite{Mun84}). Our aim is the computation of
cohomology operations on any finite simplicial
set. Steenrod squares \cite{Ste47}, Steenrod reduced powers \cite{Ste52} and Adem secondary cohomology operations
\cite{Ade52} constitute important classes of cohomology operations.

In \cite{GR98,GR99a,GR02b,GR03},
the classical definition of Steenrod cohomology operations are reinterpreted
in terms of permutations and contractions (special type of homotopy equivalences) of explicit chains for obtaining
explicit formulae
for these operations at the cochain level.
The underlying combinatorial structures of these formulae
are also studied in detail.
In \cite{GRb},
this point of view is extended to a general theory for obtaining explicit formulae
for any cohomology operation that can be
expressed at the cochain level. \cite{GRa} is devoted to give general techniques for
simplifying the underlying
simplicial structures of the formulae obtained using the process explained in
\cite{GRb}.

In this paper, we deal with the problem of the computation of Adem cohomology operations.
For this task,  we apply the machinery developed in \cite{GRb} to this special case in order
to obtain explicit formulae at the cochain level of these operations. Some properties
are also given in order to obtain simplicial ``economical'' formulae in the sense that they are
expressed only in terms of face operators. Finally, since an explicit contraction of
the cochain complex
to the cohomology can be constructed \cite{GR03}, an algorithm for computing Adem cohomology
operations at the cohomology level for nay finite simplicial set is designed in detail.

The organization of the paper follows. In Section 2, we
introduce the theoretical background from Algebraic Topology on
the concepts we use here. In Section 3,
we give a
procedure for
obtaining explicit formulae for the morphisms ${\mathcal E}_{3i+3}$ that are involved
in the definition of Adem cohomology operations  at the cochain level.
In Section 4, we design an algorithm for computing Adem cohomology operations
using the explicit formulae obtained before.

 \section{Background}\label{back}
We introduce the basic notation and terminology  we use
throughout this paper. References for the
material in this section are \cite{May67,McL75}.

A {\em simplicial set} $K$ is a graded set indexed by the
non--negative integers
 together with {\em face} and {\em degeneracy operators}
$\partial_i: K_q\rightarrow K_{q-1}$ and $s_i: K_q\rightarrow
K_{q+1}$, $0\leq i\leq q$,  satisfying the following identities:
\begin{enumerate}
\item[] $\pa_i \pa_j = \pa_{j-1} \pa_i
\mbox{ if } i<j$,
\item[]  $s_i s_j= s_{j+1} s_i
\mbox{ if } i \leq j$,
\item[]$\pa_i s_j = s_{j-1} \pa_i
\mbox{ if }i<j$, $\;\;\;\pa_j s_j = 1_{\scst K} =\partial_{j+1} s_j$,
$\;\;\;\pa_i s_j = s_j \pa_{i-1} \mbox{ if } i>j+1$.
\end{enumerate}
  The elements of $K_q$ are called $q$--{\em simplices}.
  A simplex $x$ is
  {\em degenerate} if $x=s_i y$ for some simplex $y$ and
  degeneracy operator $s_i$; otherwise, $x$ is {\em non--degenerate}.

 The {\em cartesian product} $K\ti L$ is a simplicial set whose
simplices and face and degeneracy operators are given by $(K\ti
L)_q=K_q\ti L_q$, $\pa_i(x,y)=(\pa_ix,\pa_i y)$ and $
s_i(x,y)=(s_i x,s_iy)$. The {\em semi--direct product} $
{\mathcal G}\ti_{\chi}  {\mathcal G}'$ is the group ${\mathcal
G}\ti {\mathcal G}'$ with the group--operation
$(g_1,g'_1)\cdot(g_2,g'_2)=(g_1\cdot\chi(g'_1,g_2),g'_1\cdot g'_2)\,,$
where $\chi: {\mathcal G}'\ot {\mathcal G}\ra {\mathcal G}$
satisfies that $\chi (\bar{0},g_1)=g_1,$
$\chi(g'_1,\chi(g'_2,g_1))=\chi(g'_1\cdot g'_2,g_1)$ and
$\chi(g'_1,g_1\cdot g_2)=\chi(g'_1,g_1)\cdot \chi(g'_1,g_2)$.

Let $R$ be a commutative ring with identity $1\neq 0$.
 A {\em chain complex} is a graded $R$--module $C=\oplus_{n\in{\bf Z}}\, C_n$
 together
 with an $R$--module endomorphism of degree $-1$,
 $d=\sum_{n\in {\bf Z}} d_n:C_n\rightarrow C_{n-1}$,
 such that
 $dd=d^2$ is zero. The map $d$ is called the {\em differential} of $C$.
 Ker $d_{n}$ is the module of $n$--{\em cycles} in $C$;
Im $d_{n+1}$ is the module of $n$--{\em
boundaries} in $C$; the quotient $H_{n}(C)=\mbox{Ker }
d_{n}/\mbox{Im }d_{n+1}$
is the $n^{th}$ {\em homology module} of $C$.
The homology class of a cycle $a\in$ Ker $d_n$ is denoted by
$[a]$. The $n^{th}$ homology of $C$ with coefficients in a ring
$G$ is defined by $H_n(C;G)=H_n(C\ot G)$.
Whenever two graded objects $x$ and $y$ of degree $p$ and $q$ are
interchanged we apply {\em Koszul convention} and introduce
the sign $\left( -1\right) ^{pq}.$ The {\em tensor product} of
chain complexes $C$ and $D$ is the chain complex
$C\ot D$ with differential
$d_{\scst C\ot D}=d_{\scst C}\ot 1_{\scst D}+1_{\scst C}\ot
d_{\scst D}$. Thus if $x\in C_{q}$ and $y\in D_{r},$ an
application of  Koszul convention gives
$d_{\scst C\ot D}( x\otimes y)
=d_{\scst C}( x) \ot y+( -1)^{q}x\ot d_{\scst D} (y)$.
A module homomorphism $f:C\rightarrow D$ of
degree zero such that $df=fd$ is a {\em chain map}. If
$f:C\rightarrow D$ and $g:C'\rightarrow D'$ are chain maps, so is $f\ot
g:C\ot C'\rightarrow D\ot D'$.

Let $C$ be a chain complex and $G$ an $R$--module. Form the
abelian group $C^n=\mbox{Hom}_{R}(C_n,G)$, for all $n$; its
elements are the module homomorphisms $c:C_n\rightarrow G$,
called $n$--{\em cochains} of $C^*$. The differential $d:C\ra
C$ induces a map  of degree $+1$, $\delta:C^*\ra C^*$, defined by $\delta^n
c=(-1)^{n+1}cd_{n+1}: C_{n+1}\ra G$, for all $c\in C^n$ and  for
all $n$. The {\em cohomology} of $C$ is the family of abelian
groups $H^n(C,G)=$Ker $\delta^n/$Im $\delta^{n-1}$. An element of Im $\delta^{n-1}$ is called an $n$--coboundary and an
element of Ker $\delta^n$ an $n$--cocycle.

The chain complex of a simplicial set $K$ with coefficients
in $R$, denoted by $C(K)$ is
 constructed as follows. Let $C_n(K)$ denote the free $R$--module on the set $K_n$.
The face operators $\pa_i$ yield module maps $C_n(K)\ra
C_{n-1}(K)$, which we also call $\pa_i$; their alternating sum
$d=\sum_{i\geq 0} (-1)^i\pa_i$ is the differential of $C(K)$. The
normalized chain complex $C^{\scst N}(K)$ is the chain complex
defined as the quotient
$C^{\scst N}_n(K) = C_n(K)/s(C_{n-1}(K))$,
where $s(C_{n-1}(K))$ denotes the free $R$--module on the set of
all degenerate $n$--simplices of $K$. Since we will always work
with normalized chain complexes, we simplify notation and write
$C(K)$ instead of $C^{\scst N}(K)$. The (co)homology
of $K$
  is, by definition, the (co)homology of $C(K)$.
   The
cohomology of $K$ is an algebra with the {\em cup
product}\label{cupproduct} $\smile:$
$H^i(K;G)\ot H^{j}(K;G)$ $\ra
H^{i+j}(K;G)$ defined by $[c]\smile[c']=[c\smile c']$, where $G$ is a group,
$c\in$ Ker $\delta^i$, $c'\in$ Ker $\delta^j$ and $c\smile
c'(x)=\mu(c(\pa_{i+1}\cdots\pa_{i+j}x)\ot
c'(\pa_0\cdots\pa_{i-1}x))\,,$ with $x\in C_{i+j}(K)$ and $\mu$
being the operation on $G$.
The following chain maps are needed to be defined:
\begin{itemize}
\item  The {\em diagonal map}
$\Delta:\; C(K)\ra C(K^{\ti n})$ is defined by
$\Delta\,x=(x,\stackrel{n\mbox{ \scsz times}}\dots,x)$.
\item The {\em cyclic permutations}

$t: C(K^{n})\ra C(K^{\ti n})$ such that
$t(x_1, x_2, \dots, x_n)=(x_2, \dots, x_n, x_1)$;

and $T:
C(K)^{\ot n}\ra C(K)^{\ot n}$ defined by

$T(x_1\ot x_2\ot
\cdots \ot x_n)$ $=(-1)^{|x_1|(|x_2|+ \cdots +|x_n|)}$ $ x_2\ot \cdots
\ot x_n\ot x_1$.
\end{itemize}

A {\em differential graded module} (DG--module, for short) $M$ is
a chain complex such that $M_n=0$ for all $n<0$. A {\em
DGA--module} $(M,\xi,\eta)$ (we will write it simply $M$ when no
confusion can arise) is a DG--module $M$ endowed with two
morphisms called the {\em augmentation} $\xi:\;M_0\ra R$ and the
{\em coaugmentation}, $\eta:\;R\ra M_0$. It is required that
$\xi\,\eta=1_{\scst R}$ and $\xi\,d=0$.

A {\em DGA--algebra} $(A,\mu)$ (resp. {\em DGA--coalgebra}
$(B,\nabla)$) is a DGA--module endowed with a morphism
$\mu:\,A\ot A\ra A$, called {\em product} on $A$, such that
$\mu(\mu\ot 1_{\scst A})=\mu(1_{\scst A}\ot\mu)$ and
$ \mu(\eta_{\scst A}\ot 1_{\scst A})= 1_{\scst
A}=\mu(1_{\scst A}\ot \eta_{\scst A})$. Resp. $\nabla:\, B \ra B
\ot B$, called {\em coproduct} on $B$, where $(\nabla\ot
1_{\scst B})\nabla=(1_{\scst B}\ot\nabla) \nabla$ and
$(\eta_{\scst B}\ot 1_{\scst B})\nabla= 1_{\scst
B}=(1_{\scst B}\ot \eta_{\scst B})\nabla$.

The free $R$--algebra generated by a group ${\mathcal G}$ is a
DGA--algebra denoted by
$(R[{\mathcal G}],\xi_{\scst \mathcal G},\eta_{\scst \mathcal G},\mu_{\scst \mathcal G})$
such that it is zero in each degree
except for degree zero, where $R_0[{\mathcal
G}]=\{\sum_{a\in A}\lambda_{a} a: \;\lambda_{a}\in
R\;\mbox{ and }\;A\mbox{ is a finite subset of }{\mathcal
G}\}$.
The product $\mu_{\scst \mathcal G}$, the augmentation $\xi_{\scst \mathcal G}$ and the
coaugmentation $\eta_{\scst \mathcal G}$ are given by
$\mu_{\scst
\mathcal G}((\sum \lambda_{a} a)\ot (\sum
\lambda_{a'} a')) =\sum
\lambda_{a}\lambda_{a'}(a+a')$, $\xi_{\scst \mathcal
G}(\sum \lambda_{a} a)=\sum \lambda_{a}$ and
$\eta_{\scst \mathcal G}(\lambda)=\lambda\,\0$, where
$a,a'\in \mathcal G$ and  $\;\lambda_{a}, \lambda_{a'},\lambda\in
R$.

The {\em reduced bar construction}
$\bar{B}({\mathcal G})$
 of the DGA--algebra $R[{\mathcal G}]$ is defined
  (as a graded module)
  in degree $n>0$ by $\bar{B}_n({\mathcal G})$ $=R[{\mathcal
G}]^{\ot n}/s(R[{\mathcal G}]^{\ot n})$ where $s(R[{\mathcal
G}]^{\ot n})$ is the $R$--module generated by all the elements of
$R[{\mathcal G}]^{\ot n}$ of the form $a_1\ot \cdots\ot
\0\ot\cdots\ot a_n$; and $\bar{B}_0({\mathcal G})=R$. The
element of $\bar{B}_0({\mathcal G})$ corresponding to the identity
in $R$ is denoted by $[\;]$ and an element $a_1\ot\cdots\ot a_n$
of $\bar{B}_n({\mathcal G})$ is denoted by $[a_1|\cdots|a_n]$.
 The differential of $\bar{B}({\mathcal G})$ is given by
 \begin{itemize}
 \item[] $d([a_1|\cdots|a_n])=\D\xi_{\scst \mathcal
G}(a_1)[a_2|\cdots|a_n]+ (-1)^n[a_1|\cdots|a_{n-1}]\xi_{\scst \mathcal G}(a_n)$
\\
\mbox{ }\hspace{3cm}$\D+\sum_{1\leq i\leq n-1}(-1)^i
 [a_1|\cdots|a_{i-1}|\mu_{\scst \mathcal G}(a_i\ot a_{i+1})|a_{i+2}|\cdots|a_n]\,.$
 \end{itemize}
 Observe that $d[a_1]=0$ for all $[a_1]\in \bar{B}_1({\mathcal G})$.
The augmentation and the coaugmentation on $\bar{B}({\mathcal
G})$ coincide with the
 identity on $R$. Moreover,
$\bar{B}({\mathcal G})$ is a DGA--coalgebra
 with the coproduct: $$ \nabla([a_1|\cdots|a_n])=\sum_{0\leq i\leq
 n}[a_1|\cdots|a_i]\ot [a_{i+1}|\cdots|a_n].$$
Let $(B,\nabla)$ be a DGA--coalgebra and  $(A,\mu)$ a
DGA--algebra. A {\em twisting cochain}
$\kappa$, is a graded module morphism of degree $-1$, $\kappa:\,B\ra A$,
 satisfying that
$d_{\scst A}\kappa+\kappa\,d_{\scst B}+\mu(\kappa\ot
\kappa)\nabla=0$, $\xi_{\scst A}\kappa=0$ and $\kappa\,\eta_{\scst B}=0$.
Let  $M$ be an $A$--DG--module
(where $\nu:M\ot A\ra M$ is the (right) $A$--module structure on
$A$). Define the morphism $d_{\kappa}:\, M\ot B\ra M\ot
B$ where $d_{\kappa}(m\ot b) =(d_{\scst M\ot B}+\kappa\cap)(m\ot
b)$ and $\kappa\cap=(\nu\ot 1_{\scst B})(1_{\scst M}\ot
\kappa\ot 1_{\scst B}) (1_{\scst M}\ot \nabla)$. The graded module
$M\ot B $ endowed with $d_{\kappa}$ is a DG--module
denoted by $M\ot_{\kappa} B$ and called {\em twisted tensor
product} by the twisting cochain  $\kappa$. An example of twisted
tensor product is $ R[{\mathcal
G}]\ot_{\theta}\bar{B}({\mathcal G})$, where the twisting
cochain $\theta$ called the {\em universal twisting cochain} is
given by $\theta([a_1])=a_1-\0$ and $\theta([a_1|\cdots|a_n])=0$ for $n>1$.

We deal with an special type of homotopy equivalences. A {\em
contraction} $r$ of a DG--module $N$ to a DG--module $M$,
consists in three morphisms $(f,g,\phi)$ where
 $f:\,N\ra M$ (projection) and $g:\,M\ra N$ (inclusion) are DG-module morphisms of degree zero,
and $\phi:\,N\ra N$ (homotopy operator) is a morphism of degree $-1$ satisfying that
$fg= 1_{\scst M}$, $\phi d + d\phi = g f-1_{\scst N}$. Moreover, it is  required that
$\phi g =0$, $f\phi = 0$, $\phi\phi = 0$.
A contraction will be denoted by $r=(f,g,\phi):\,N\Ra M$ or
briefly $N\stackrel{\scst r}\Ra M$. Note that the importance of
having this structure of $M$ to $N$ is that $N$ is
``smaller'' than $M$ although both have the same
homology. Let $r=(f,\, g,\,\phi):\,N\Ra M$ and $r'=(f',\,
g',\,\phi'):\,N'\Ra M'$
  be two contractions,  then  the following
contractions can be constructed:
\begin{itemize}
\item The tensor product contraction:
$r\ot r'=(f\ot f', g\ot g', \phi\ot g' f'+1_{\scst N}\ot
\phi'): N\ot N'\Ra M\ot M'$.
\item If $M=N'$, the  composition contraction:
$r' r=(f' f, g g', \phi+g\phi' f): N\Ra  M'$.
\end{itemize}

Let $p$ and $q$ be  non--negative integers, a {\em
    $(p,q)$--shuffle} $(\al, \be)$ is a partition of the set
    $\{0,1, \ldots, p+q-1\}$ in two disjoint
    subsets, $\al_1 < \cdots < \al_p$ and
    $\be_1 < \cdots < \be_q$, of $p $ and $q$ integers,
    respectively. The signature  of the
    shuffle $(\al, \be)$ is defined by
$sig(\al, \be)= \sum_{1\leq i\leq p} \al_i - (i-1)$.

  An Eilenberg--Zilber contraction \cite{EZ59}  of
$C(K\ti L)$ to $C(K)\ot C(L)$, where $K$ and $L$ are
simplicial sets, is a triple $r_{\scst EZ}=(Aw,Em,Sh)$ where:
\begin{itemize}
\item  Alexander--Whitney operator
$Aw:\,C(K \ti L) \ra C(K)\ot C(L)$ is defined by:
$\D Aw(x_m
\ti y_m)= \sum _{0\leq i\leq m} \pa_{i+1} \cdots \pa_m x_m \ot
\pa_0 \cdots \pa_{i-1} y_m\,.$

\item Eilenberg--Mac Lane operator
$Em:\, C(K)\ot C(L)\ra C(K \ti L)$ by:

$\D Em(x_p \ot y_q) =
\sum_{( \al, \be ) \in \{ (p,q)- \mbox{\scriptsize shuffles}
\}} (-1)^{sig( \al , \be )}(s_{\be} x_p,\, s_{\al} y_q)$\\
where
$s_{\be}$ $=s_{\be_q}\cdots s_{\be_1}$ and $s_{\al}=s_{\al_p}\cdots
s_{\al_1}$.
\item Shih operator
$Sh:\,C(K \ti L)\ra C(K \ti L)$ by:

$\D Sh((x_m,\, y_m)) = \sum_{T(m)}(-1)^{sg}(s_{\bb+\bar{m}} \pa_{m-q+1} \cdots \pa_m x_m,
s_{\al+\bar{m}}\pa_{\bar{m}} \cdots \pa_{m-q-1}y_m)$
where $sg= \bar{m}+sig( \al , \be)+1$,
$\bar{m}=m-p-q$,
$\bb+\bar{m}$ $=\{\be_q+\bar{m},\cdots\be_1+\bar{m},\bar{m}-1\}$,
$\al+\bar{m}$ $=\{\al_{p+1}+\bar{m},\cdots\al_1+\bar{m}\}$
 and $T(m)$ $=\{0\leq p\leq m-q-1\leq
m-1,\;(\al,\be)\in\{(p+1,q)\mbox{--shuffles}\}\}$.
\end{itemize}
 It is evident that  Alexander--Whitney operator has a
polynomial nature.
However, Eilenberg--MacLane and Shih operator have an essential
 ``exponential'' character because
 shuffles  of degeneracy operators are involved in their respective formulation.
 A recursive
 formula for the Shih operator is given in \cite{EM54}.  An explicit formula
 for this operator was stated by J. Rubio
 \cite{Rub91} and proved by F. Morace \cite{Rea00}.

It is possible to construct  a contraction $r_{\scst
EZ(p)}=(Aw_p,Em_p,Sh_p)$
of $ C (K^{\ti p})
$ to  $C(K)^{\ot p}$ ($p\geq 0$),
appropriately composing  Eilenberg--Zilber contractions.

Note that the cup product defined on page \pageref{cupproduct}
can be written as $c\smile c'=\mu (c\ot c')Aw\Delta$ where
$\mu$ is the operation defined on $G$.
Generalizations of the cup product are the cochain
mappings called
 {\em cup--$i$ product} \cite{Ste47}, $\smile_i:$ $C^n(K;\Z)$
 $\ot C^m(K;{\bf Z})\ra C^{n+m-i}(K;\Z)$ given by
$c\smile_i c'=\mu (c\ot c')D_i\Delta$ where $D_i=Aw (t\,Sh)^i$.
The following relation holds (up to sign): \begin{equation}
\label{pie} \delta (c\smile_i c')=\pm c\smile_{i-1}c'\pm
c'\smile_{i-1}c +\delta c\smile_i c'+c\smile_i\delta
c'\,.\end{equation} Taking $[c]\in H^j(K;\Z)$ and defining
$Sq^i[c]=[c\smile_{j-i} c]$, then $Sq^i(c)$ $\in$ $
H^{j-i}(K;\Zr)$. These cohomology operations are called {\em
Steenrod squares} \cite{Ste47}.

Explicit formulae for $D_i$ and therefore for computing cup--$i$ products,
Steenrod squares and Steenrod reduced powers \cite{Ste52}
(a generalization of
Steenrod squares), can be found in
\cite{GR98,GR99a,GR99b,GR02b}.

In \cite{Ade52}, J. Adem constructed secondary cohomology
operations using  relations on iterated Steenrod squares. He
proved the relation $Sq^2Sq^2+Sq^3Sq^1=0$ by means of the existence of
cochain mappings ${\mathcal E}_j:C^p(K^{\ti 4};\Z)\ra
C^{p-j}(K;\Zr)$
 such that mod $2$
\begin{equation}\label{relacion1}
(c\smile_i c)\smile_{i+2}(c\smile_i c)+ (c\smile_{i+1}
c)\smile_i(c\smile_{i+1} c)=\delta {\mathcal E}_{3i+3}(c^4)\end{equation}
where $\smile_k$ is the cup--$k$ product, $c$ is a
$q$--cocycle and $i=q-2$. If $c$ is a $q$--cocycle
such that $Sq^2 (c)$ is a coboundary (that is, there exists a
cochain $b$ such that $c\smile_i c=\delta b$), then
$$ w=b\smile_{i+1} b+
b\smile_{i+2}\delta b+ {\mathcal
E}_{3i+3}(c^4)+\eta(c)\smile_{i-1}\eta(c) +\eta(c)\smile_i \delta
\eta (c)$$ is a mod $2$ cocycle, where $
\eta(c)=\frac{1}{2}(c\smile_{i+2} c)+c$.
Adem secondary cohomology operations are defined as
$\Psi_q[c]=[w]+Sq^2H^{q+1}(K;\Z)\in H^{q+3}(K;\Zr)$.

In Section 4 we give a procedure for
computing the Adem secondary cohomology operation $\Psi_q (\alpha)$ for any  integer $q$
and any cohomology class $\alpha$. For doing this,
the obtention of explicit formuale for ${\mathcal E}_{3i+3}$ is essential. We will obtain these
formulae in the following section.

\section{Adem  Cocyclic Operations}\label{secadem}

In order to obtain ``economical'' formulae for ${\mathcal E}_{3i+3}$
in the sense that they are
written only in terms of face operators, we first
write ${\mathcal E}_{3i+3}$ in terms of the component morphisms of Eilenberg--Zilber
contractions. After this, we make a simplification of the given formulae based on the fact
that
any composition of face and degeneracy operators can be put in a
unique form.

For the first aim,
we will use the following result whose demonstration
is given in \cite{GRb}.

\begin{lem}\label{lema1}
Let ${\mathcal G}$ be a group. Let $M$ and $N$ be two
$R[{\mathcal G}]$--DG--modules, where $\nu: M\ot R[{\mathcal
G}]\ra M$ and $\nu': N\ot R[{\mathcal G}]\ra N$ are the (right)
$R[{\mathcal G}]$--module structures on $R[{\mathcal G}]$.
Let  $r=(f,g,\phi):\;M\Ra N$ be a
contraction such that $(g\ot 1_{\scst
\bar{B}[{\mathcal G}]})\nu'=\nu g$.

Then $
( r\ot 1_{\scst \bar{B}[{\mathcal G}]})_{\scst \tc}
=((f\ot 1_{\scst \bar{B}[{\mathcal G}]})_{\scst \tc},\,g\ot 1_{\scst \bar{B}[{\mathcal G}]},\,
(\phi\ot 1_{\scst \bar{B}[{\mathcal G}]})_{\scst \tc})$ of
$M\ot_{\scst \theta}\bar{B}({\mathcal G})$ to $N\ot_{\scst \theta}\bar{B}({\mathcal G})
$ is a new contraction
where $ (f\ot 1_{\scst \bar{B}[{\mathcal G}]})_{\scst \tc}=\sum_{i\geq 0}
(f\ot 1_{\scst \bar{B}[{\mathcal G}]})(\tc(\phi\ot 1_{\scst \bar{B}[{\mathcal G}]}))^i$
and $ (\phi\ot 1_{\scst \bar{B}[{\mathcal G}]})_{\scst \tc}=-\sum_{i\geq 0}
((\phi\ot 1_{\scst \bar{B}[{\mathcal G}]})\tc)^i(\phi\ot 1_{\scst \bar{B}[{\mathcal G}]})$.
\end{lem}

Now, let ${\mathcal G}$ be the semi--direct product $\Zr^{\ti 2}\ti_{\chi} \Zr$
where $\chi((a,b),\1)=(b,a)$. Let $M=C(K^{\ti 4})$ and
$N=C(K)^{\ot 4}$ be two $\Zr[{\mathcal
G}]$--modules:
\begin{itemize}
\item[] $\nu(x,a_1)=z(x)=
(x_1,x_3,x_2,x_4)\,,
\;\;\;\; \nu'(y,a_1)=z'(y)=y_1\ot y_3\ot y_2\ot y_4\,,$\\
$\nu(x,a_2)=t^{\ti 2}(x) =(x_2,x_1,x_4,x_3)\,,
\;\nu'(y,a_2)=T^{\ot 2}(y)=y_2\ot y_1\ot y_4\ot y_3\,,$\\
$\nu(x,a_3)=t( x)=(x_3,x_4,x_1,x_2)\,,
\;\;\;\; \nu'(y,a_3)=T(y)=y_3\ot y_4\ot y_1\ot y_2\,;$
\end{itemize}
 where
$a_1=((\bar{0},\bar{0}),\bar{1})$,
$a_2=((\bar{1},\bar{0}),\bar{0})$,
 $a_3=((\bar{0},\bar{1}),\bar{0})$,
 $x=(x_1,x_2,x_3,x_4)\in C(K^{\ti 4})$ and $y=y_1\ot
y_2\ot y_3\ot y_4\in C(K)^{\ot 4}$. Let $r$ be the
Eilenberg--Zilber  contraction
$r_{\scst
EZ(4)}=(r_{\scst EZ}^{\ot 2})r_{\scst EZ}=(Aw_4,Em_4,Sh_4):
C(K^{\ti 4})\Ra
C(K)^{\ot 4}$ which commutes
with the structures given above. Observe that
$Aw_4=(Aw^{\ot 2})Aw$, $Em_4=Em(Em^{\ot 2})$ and $Sh_4= Sh+Em(Sh\ot Em
Aw+1\ot Sh)Aw)$.
Lemma \ref{lema1} produces the
new contraction: $$(r_{\scst EZ(4)}\ot 1_{\scst \bar{B}(\Zr^{\ti
2}\ti_{\chi}\Zr)})_{\tc}:  C(K^{\ti
4})\ot_{\theta} \bar{B}(\Zr^{\ti 2}\ti_{\chi}\Zr)\Ra C(K)^{\ot
4}\ot_{\theta}  \bar{B}(\Zr^{\ti 2}\ti_{\chi}\Zr)\,.$$
As we will see now, Adem relation (\ref{relacion1}) can be obtained from the fact
that the projection $(Aw_4\ot 1_{\scst \bar{B}(\mathcal G)})_{\tc}$ is a DG--module morphism.

\begin{thm}\label{principal}
Adem relation (\ref{relacion1}) is obtained from the identity:
\begin{itemize}
\item[] $\mu c^{\ot 4}(1_{\scst N}\ot \xi_{\scst \bar{B}({\mathcal G})}) (Aw_4\ot
1_{\scst \bar{B}({\mathcal G})})_{\tc}(d_{\scst M\ot \bar{B}({\mathcal G})}+\tc)
(\Delta(x)\ot e_{3i+3})$\\
$= \mu c^{\ot 4}(1_{\scst N}\ot \xi_{\scst \bar{B}({\mathcal G})}) (d_{\scst N\ot
\bar{B}({\mathcal G})}+\tc)(Aw_4\ot 1_{\scst \bar{B}({\mathcal G})})_{\tc} (\Delta(x)\ot
e_{3i+3})$\end{itemize}
where $\mu$ is the product on ${\bf Z}_2$,
$c\in$ Ker $\delta^q$, $x\in C_{q+3}(K)$
and  $e_{3i+3}\in \bar{B}({\mathcal G})$ is:
 \begin{itemize}
 \item[] $e_{3i+3}=e_{(3,0)}\;$ for $i=0$,
 \item[] $e_{3i+3}=e_{(3i+3,i)}+e_{(3i+3,i-1)}\;$ for $i=2,4$,
  \item[]  $e_{3i+3}=e_{(3i+3,i)}\;$
  for $i=1,3,5$ and for $i\geq 6$, $i$ being odd and $\frac{i-7}{2}$ being even,
 \item[]  $e_{3i+3}=e_{(3i+3,i)}+e_{(3i+3,i-2)}\;$
   for  $i\geq 6$, $i, \frac{i-7}{2}$ being odd,
 \item[]  $e_{3i+3}=e_{(3i+3,i)}+e_{(3i+3,i-1)}+e_{(3i+3,i-2)}\;$
  for  $i\geq 6$, $i$ being even,
 \end{itemize}
such that
\begin{itemize}
\item[] $\D e_{(3i+3,\ell)}=
\sum_{\scst S(3i+2,\ell)}[(c_{\pi_1},\0)|\dots|(c_{\pi_j},\0)|((\0,\0),\1)|(\bar{c}_{\pi_{j+1}},\0)|\dots|
(\bar{c}_{\pi_{3i+2}},\0)]$
\end{itemize}
where
$S(3i+2,\ell)=\{\pi\in \{(3i+2-\ell,\ell)\mbox{\--sh.}\}\}\cup
\{j:\,0\leq j\leq 3i+2\}$,
$[c_1,\dots,c_{3i+2}]=[(\1,\0),\stackrel{3i+2-\ell}\dots,(\1,\0),(\0,\1),\stackrel{\ell}\dots,
(\0,\1)]$ and if $c=(a,b)$ then $\bar{c}=(b,a)$.
\end{thm}
The proof of this theorem is given in an appendix at the end of this paper.
 \begin{cor}
The explicit formulae for ${\mathcal E}_{3i+3}$ are:
\begin{itemize}
\item[] ${\mathcal E}_{3}=\tilde{\mathcal E}_{(3,0)}+c\smile_0c\smile_1c\smile_2 c\;$
for $i=0$.
\item[] ${\mathcal E}_{6}=\tilde{\mathcal E}_{(6,1)}+
c\smile_1c\smile_2c\smile_3 c\;$ for $i=1$.
\item[] ${\mathcal E}_{9}=\tilde{\mathcal E}_{(9,2)}+\tilde{\mathcal E}_{(9,1)}+
c\smile_2c\smile_3c\smile_4 c+c\smile_0 c\smile_7c\smile_2\;$
for $i=2$.
\item[]  ${\mathcal E}_{12}=\tilde{\mathcal E}_{(12,3)}+
c\smile_3c\smile_4c\smile_5 c+c\smile_2c\smile_7c\smile_3 c
+c\smile_3c\smile_8c\smile_1 c$\\ \mbox{ }\hspace{0.2cm} for $i=3$.
\item[] ${\mathcal E}_{15}=\tilde{\mathcal E}_{(15,4)}
+\tilde{\mathcal E}_{(15,3)}+
c\smile_4c\smile_5c\smile_6 c\;$ for $i=4$.
\item[]  ${\mathcal E}_{18}=\tilde{\mathcal E}_{(18,5)}+
c\smile_5c\smile_6c\smile_7 c+c\smile_4c\smile_9c\smile_5 c
+c\smile_4c\smile_{11}c\smile_3 c$\\ \mbox{ }\hspace{0.2cm} for $i=5$.
\item[]
${\mathcal E}_{3i+3}=\tilde{\mathcal E}_{(3i+3,i)}+\tilde{\mathcal E}_{(3i+3,i-1)}+
\tilde{\mathcal E}_{(3i+3,i-2)}+c\smile_i c\smile_{i+1} c\smile_{i+2}$\\
\mbox{ }\hspace{0.2cm} for $i\geq 6$, $i$ being even, $\frac{i-6}{2}$ being odd.
\item[]
${\mathcal E}_{3i+3}=\tilde{\mathcal E}_{(3i+3,i)}+\tilde{\mathcal E}_{(3i+3,i-1)}+
\tilde{\mathcal E}_{(3i+3,i-2)}+c\smile_i c\smile_{i+1} c\smile_{i+2}$\\
\mbox{ }\hspace{0.5cm} $+c\smile_i c\smile_{i+5} c\smile_{i-2}\;$
for $i\geq 6$, $i,\frac{i-6}{2}$ being even.
\item[] ${\mathcal E}_{3i+3}=\tilde{\mathcal E}_{(3i+3,i)}+\tilde{\mathcal E}_{(3i+3,i-2)}+
c\smile_i c\smile_{i+1} c\smile_{i+2}+
c\smile_i c\smile_{i+4} c\smile_{i-1}$\\
\mbox{ }\hspace{0.5cm} $+c\smile_i c\smile_{i+5} c\smile_{i-2}+
c\smile_{i-1} c\smile_{i+6} c\smile_{i-2}$\\
\mbox{ }\hspace{0.2cm} for $i\geq 7$, $i,\frac{i-7}{2}$ being odd.
\item[]
${\mathcal E}_{3i+3}=\tilde{\mathcal E}_{(3i+3,i)}+
c\smile_i c\smile_{i+1} c\smile_{i+2}+c\smile_i c\smile_{i+4} c\smile_{i-1}$\\
\mbox{ }\hspace{0.5cm} $+c\smile_i c\smile_{i+5} c\smile_{i-2}
+c\smile_{i-1} c\smile_{i+6} c\smile_{i-2}$\\
\mbox{ }\hspace{0.2cm} for $i\geq 7$, $i$ being odd,  $\frac{i-7}{2}$ being even.
\end{itemize}
where $\tilde{\mathcal E}_{(3i+3,\ell)}$ denotes the composition:
\begin{equation}\label{pozi2}\begin{array}{l}
\sum \mu c^{\ot 4}(Aw_4(t\,Sh_4)^{m_0}(t^{\ti
2}Sh_4)^{n_1}(t\,Sh_4)^{m_1}\cdots z\,Sh_4\cdots (t^{\ti 2}Sh_4)^{n_j}
(t\,Sh_4)^{m_{j}}\qquad\\
\;+Aw_4(t^{\ti 2}Sh_4)^{m_0}(t\,Sh_4)^{n_1}(t^{\ti 2}Sh_4)^{m_1}
\cdots z\,Sh_4\cdots (t\,Sh_4)^{n_j}(t^{\ti 2}\,Sh_4)^{m_j}\Delta(x)
\end{array}\end{equation}
such that  the sum is taken over the set
$\{(n_1,\dots,n_j,m_0,\dots,m_{j}):$  $n_1+\cdots +n_j=3i+2-\ell$,
$m_0+\cdots+m_j=\ell$, $n_k\geq 0$, $m_k\geq 0 \}$.
\end{cor}
Observe that an algorithm designed for computing Adem operations from these formulae
 would be too
slow for practical implementation because  the morphisms $Sh_4$ and $Em_4$
 present an exponential number of summands in their formuale.
For this reason, the idea of
simplification arises in a natural way.
In the following subsection we give some useful properties  for the normalization process.
Some work have been  done in this way by the authors in \cite{GRa}.

\subsection{Adem Cocyclic Operations and Combinatorics}\label{p}

A refinement of ${\mathcal E}_{3i+3}$, only in terms of face operators
of $K$, is feasible using  a normalization process.
The following properties will be useful in this process.
\begin{prop} Let $K$ be a simplicial set then:
\begin{enumerate}
\item[1.] Any composition of face and degeneracy operators of $K$ can be put
in the unique form:
$s_{j_t}\cdots s_{j_1}\partial_{i_1}\cdots \partial_{i_s}$,
where $j_t > \cdots >j_1\geq 0$ and $i_s > \cdots > i_1 \geq 0$.
\item[2.] Those summands of the simplified formula for ${\mathcal E}_{3i+3}$ with a factor
having a degeneracy operator in its expression are null.
\item[3.]  All the summands of
$(1\ot Sh)Aw z\,Sh: C(K^{\ti 4})\ra C(K^{\ti 2})^{\ot 2}$ are null.
\item[4.] $Aw\,z\,Em=(Em\ot Em)z'(Aw\ot Aw)$ where $Aw\,z\,Em:C(K^{\ti 2})^{\ot 2}\ra C(K^{\ti 2})^{\ot 2}\,.$
\item[5.] It is satisfied that
$Sh \,t\,Sh=Sh\,t\tilde{Sh}: C(K^{\ti 2})\ra C(K^{\ti 2})$ where $\tilde{Sh}:C(K^{\ti 2})\ra
C(K^{\ti 2})$  consists in all the summands of $Sh$ with
$\be_q<\al_{1}$.
\end{enumerate}
\end{prop}
Let us prove the third assertion. Let $(x_1,x_2,y_1,y_2)\in
C_m(K^{\ti 4})$ then the second factor of each summand of $(1\ot
Sh)Aw\, z\,Sh(x_1,x_2,y_1,y_2)$ has the form
\begin{itemize}
\item[] $ (s_{\scst \bar{b}+\bar{n}}\pa_{\scst n-q'+1}\cdots\pa_{\scst n}
\pa_{\scst 0}\cdots\pa_{\scst i-1}s_{\scst \bb+\bar{m}}\pa_{\scst
m-q+1}
\cdots\pa_{\scst m} x_{\scst 2},$\\
$s_{\scst a+\bar{n}}\pa_{\scst \bar{n}}\cdots \pa_{\scst n-q'-1}
\pa_{\scst 0}\cdots\pa_{\scst i-1}s_{\scst \al+\bar{m}}\pa_{\scst
\bar{m}}\cdots \pa_{\scst m-q-1}y_{\scst 2})$\end{itemize}
where $n=m-i+1$, $0\leq i\leq m+1$, $0\leq p\leq
m-q-1\leq m-1$, $0\leq p'\leq n-q'-1\leq n-1$,
$(\al,\be)\in \{(p+1,q)\mbox{--sh.}\}$ and $(a,b)\in \{(p'+1,q')\mbox{--sh.}\}$.

If $i=m+1$ then $n=0$,
and  $Sh$ is $0$ in degree $0$. When $0\leq i\leq m$ then
$\al+\bar{m}\,\cap\, \{m-q'+1,m-q'+2\dots,m\}=\emptyset$ because in other case the
corresponding cartesian product is degenerate. Then
$\{m-q'+1,m-q'+2\dots,m\}\subset \bb+\bar{m}$. For the same reason than above,
 $\bb+\bar{m}\,\cap\, \{m-p'-q',m-p'-q'+1,\dots,m-q'\}=\emptyset$, and it follows that
 $\{m-p'-q',m-p'-q'+1,\dots,m-q'\}\subset   \al+\bar{m}$ and $\bar{m}-1<m-p'-q'$.
But this case is also degenerate.

We leave it  to the reader to verify the other
assertions because they follow the same technique.

For example, the normalized formula for ${\mathcal E_3}$ is:
\begin{itemize}
\item[] $
{\mathcal E_3}(c^4)(x)=\mu c^{\ot 4}
(\pa_1\pa_4\pa_5\ot \pa_3\pa_4\pa_5\ot \pa_0\pa_1\pa_2\ot\pa_0\pa_1\pa_4\\
+\pa_1\pa_2\pa_3\ot \pa_0\pa_1\pa_2\ot \pa_3\pa_4\pa_5\ot \pa_3\pa_4\pa_5
+\pa_2\pa_3\pa_4\ot\pa_0\pa_1\pa_2\ot \pa_0\pa_4\pa_5\ot\pa_0\pa_4\pa_5\\
+\pa_3\pa_4\pa_5\ot \pa_0\pa_1\pa_3\ot \pa_0\pa_1\pa_5\ot \pa_0\pa_1\pa_5
+\pa_3\pa_4\pa_5\ot \pa_0\pa_1\pa_4\ot \pa_0\pa_1\pa_2\ot \pa_0\pa_1\pa_2)\Delta (x)
$
\end{itemize}
where $c$ is a $2$--cocycle,  $x$ is a $5$--simplex
and $\mu$ is the product on ${\bf Z}_2$.

\section{Adem Cohomology Operations}\label{homo}

An algorithm  for computing the homology (if it is torsion free or the
ground ring is a field)
of a finite simplicial set $K$ and a contraction $(f,g,\phi)$
of $C(K)$
to $H(K)$  appear in
\cite{GR03}.
 The complexity of the algorithm  is $O(m^3)$, where $m$ is the number of simplices of
$K$.
From this contraction, it is easy to derive another one, $(f^*,g^*,\phi^*)$, of
$C^*(K;{\bf Z})$ to $H^*(K;{\bf Z})$ (see \cite{GR03}).
An interesting property of this contraction is that
if $c$ is a coboundary, then $c=\delta\phi^*(c)$.

For attacking the computation of Adem secondary cohomology operations,
we will see  that the homotopy operator $\phi^*$,
the explicit formulae for computing cup--$i$ products given in \cite{GR99a} and the
explicit formulae for ${\mathcal E}_{3i+3}$
 are
essential. The steps for computing $\Psi_q$ are the following.
\begin{alg} {\em Algorithm for computing Adem cohomology operations.}

 {\sc Input:} A finite simplicial set $K$ and an integer $q=i-2\geq 2$.
\begin{enumerate}
\item[1.] Compute $(f^*,g^*,\phi^*)$ of $C^*(K;{\bf Z})$ to $H^*(K;{\bf Z})$.\\
Let $\{\alpha_1,\dots,\alpha_p\}$ be a set of generators of
$H^q(K;{\bf Z})$.
\item[2.] Compute $Sq^2: H^q(K;{\bf Z})\ra H^{q+2}(K;{\bf Z}_2)$ as follows:\\
For $j=1$ to $j=p$, compute $f^*(g^*(\alpha_j)\smile_i g^*(\alpha_j))$.
\item[3.]Put the matrix corresponding to $Sq^2: H^q(K;{\bf Z})\ra H^{q+2}(K;{\bf Z}_2)$
in a digonal form $D$.\\
Let $\{\beta_1,\dots,\beta_r \}$ be a set of generators of Ker
$Sq^2$ obtained using $D$.
\item[4.] For $j=1$ to $j=r$, compute $c_j=g^*(\beta_j)$, $b_j=\phi^*(c_j\smile_i c_j)$ and
$\eta(c_j)=\frac{1}{2}(c_j\smile_{i+2} c_j+c_j)$.
\item[5.] For $j=1$ to $j=r$, compute $w_j=b_j\smile_{i+1} b_j+
b_j\smile_{i+2}\delta(b_j)+{\mathcal E}_{3i+3}(c_j^4)+\eta(c_j)\smile_{i-1}\eta(c_j)+
\eta(c_j)\smile_i\delta\eta(c_j)$.
\item[6.] For $j=1$ to $j=r$, compute $f^*(w_j)$.
\end{enumerate}
{\sc Output:} $\Psi_q: H^q(K;{\bf Z})\ra H^{q+3}(K;{\bf Z}_2)$
\end{alg}
A practical implementation in $\mbox{\em Mathematica}^{\mbox{\scsz \copyright}}$ and a  concrete example about the computation of
the first Adem operation
$\Psi_2$ is given in  \cite{GR02a}.

\bibliographystyle{elsart-harv}

\section*{Appendix: Proof of Theorem \ref{principal}}

First of all, the chain $e_{3i+3}$ ($i\geq 0$) is obtained using
 the following composition of contractions:
$$\begin{array}{r}
\bar{B}(\Zr\ti_{\chi}\Zr^{\ti 2})
\stackrel{\scst r_{\scst EZ(\chi)}}\Ra
\bar{B}(\Zr)\ot_{\theta}\bar{B}(\Zr^{\ti 2})
\stackrel{\scst (1\ot r_{\scst \bar{B}})_{\scst \tc}}\Ra\bar{B}(\Zr)\ot_{\theta}\bar{B}(\Zr)^{\ot 2}\\
\Downarrow{\scst \scst r_{\scst 2}^{\scst \ot 3}}\\
E(u)\ot \Gamma(v)\ot_{\theta}(E(u_1)\ot\Gamma(v_1)\ot E(u_2)\ot\Gamma(v_2))
\end{array}
$$
 where
$r_{\scst EZ( \chi)}=(Aw_{\chi}, Em_{\chi},Sh_{\chi})$  is a perturbed
Eilenberg--Zilber contraction
such that
 the inclusion $Em_{\chi}$ is not perturbed
\cite{Arm99};
$(1\ot r_{\scst \bar{B}})_{\scst \tc}$
satisfies that the inclusion is $1\ot g_{\scst \bar{B}}$ (Lemma \ref{lema1}) and
$r_2$ is an isomorphism with explicit formulae given in \cite{EM52}.

Denote the composition $Em_{\chi} (1\ot
g_{\scst\bar{B}})g_2^{\ot 3}$ of $E(u)\ot
\Gamma(v)\ot_{\theta}(E(u_1)\ot\Gamma(v_1)\ot
E(u_2)\ot\Gamma(v_2))$ to $\bar{B}(\Zr\ti_{\chi}\Zr^{\ti 2})$  by
$g_{\scst E\Gamma}$. Let $w_i^j$ (i=1,2) be the generator of
$E(u_i)\ot\Gamma(v_i)$ in degree $j$. Then,
\begin{itemize}
\item[] $e_{(3i+3,\ell)}=g_{\scst E\Gamma}(u\ot
w_1^{\ell}\ot w_2^{3i+2-\ell})$\\
\mbox{ }\hspace{1.2cm}$=\D\sum_{\scst S(3i+2,\ell)}[(c_{\pi_1},\0)|\dots|(c_{\pi_j},\0)|((\0,\0),\1)|(\bar{c}_{\pi_{j+1}},\0)|\dots|
(\bar{c}_{\pi_{3i+2}},\0)]$.
\end{itemize}

We have to prove that Adem relation (\ref{relacion1}) is obtained from the identity:
\begin{equation}
\label{fact1}\begin{array}{c}
\mu c^{\ot 4}(1_{\scst N}\ot \xi_{\scst \bar{B}({\mathcal G})}) (Aw_4\ot
1_{\scst \bar{B}({\mathcal G})})_{\tc}(d_{\scst M\ot \bar{B}({\mathcal G})}+\tc)
(\Delta(x)\ot e_{3i+3})\\
= \mu c^{\ot 4}(1_{\scst N}\ot \xi_{\scst \bar{B}({\mathcal G})}) (d_{\scst N\ot
\bar{B}({\mathcal G})}+\tc)(Aw_4\ot 1_{\scst \bar{B}({\mathcal G})})_{\tc} (\Delta(x)\ot
e_{3i+3})\end{array}\end{equation}

Since $e_{3i+3}$ is defined as a sum of elements of the form
$e_{(3i+3,\ell)}$ for $\ell=i, i-1,$ $i-2$,
the left--hand side of (\ref{fact1}) has the form:
\begin{equation}
\label{one}\sum_{\ell}\mu c^{\ot 4}(1_{\scst N}\ot \xi_{\scst \bar{B}({\mathcal G})}) (Aw_4\ot
1_{\scst \bar{B}({\mathcal G})})_{\tc}(d_{\scst M}\ot 1_{\scst \bar{B}({\mathcal G})})
(\Delta(x)\ot e_{(3i+3,\ell)})\end{equation}
and the right--hand side,  the form:
\begin{equation}
\label{two}
\sum_{\ell}\mu c^{\ot 4}(1_{\scst N}\ot \xi_{\bar{B}({\mathcal G})}) (Aw_4\ot
1_{\scst \bar{B}({\mathcal G})})_{\tc} (\Delta(x)\ot b_{(3i+2,\ell)})\end{equation}
where
\begin{itemize}
\item[]$b_{(3i+2,\ell)}=
g_{\scst E\Gamma}d_{\scst \bar{B}({\bf Z}_2\ti_{\chi}{\bf Z}_2^{\ti 2})}(e_{(3i+3,\ell)})=
g_{\scst E\Gamma}(w_1^{\ell}\ot w_2^{3i+2-\ell}+
w_1^{3i+2-\ell}\ot w_2^{\ell})$\\
$\D =\sum_{S(3i+2,\ell)}
[(c_{\pi_1},\0)|\dots|(c_{\pi_{3i+2}},\0)]+
[(\bar{c}_{\pi_1},\0)|\dots| (\bar{c}_{\pi_{3i+2}},\0)]$
\end{itemize}
On one hand, (\ref{one}) can be simplified to
$\delta(\sum_{\ell} \tilde{\mathcal E}_{(3i+3,\ell)})$
where $\tilde{\mathcal E}_{(3i+3,\ell)}$ denotes the composition (\ref{pozi2}) defined
on page \pageref{pozi2}.

On the other hand, (\ref{two}) can be written mod $2$ as
\begin{equation}\label{pono}\begin{array}{l}
\D \sum_{\ell}\sum \mu c^{\ot 4}(Aw_4(t\,Sh_4)^{m_0}(t^{\ti
2}Sh_4)^{n_1}(t\,Sh_4)^{m_1}\cdots (t^{\ti 2}Sh_4)^{n_j}
(t\,Sh_4)^{m_{j}}\\
\mbox{ }\qquad+Aw_4(t^{\ti 2}Sh_4)^{m_0}(t\,Sh_4)^{n_1}(t^{\ti 2}Sh_4)^{m_1}
\cdots (t\,Sh_4)^{n_j}(t^{\ti 2}\,Sh_4)^{m_j})\Delta(x)
\end{array}\end{equation}
where the last sum is taken over the set
$\{(n_1,\dots,n_j,m_0,\dots,m_{j}):$  $n_1+\cdots +n_j=3i+2-\ell$,
$m_0+\cdots+m_j=\ell$, $n_k\geq 0$, $m_k\geq 0 \}$.

In order to simplify (\ref{pono}),
we will use the facts that

 $Aw_4 (t\,Sh_4)^n=D_n\;$ if $n\geq 0$.

 $Aw_4 (t\, Sh_4)^n t^{\ti 2}Sh_4=0\;$ if $n\geq 1$.

 $\D Aw_4 (t^{\ti 2}Sh_4)^n
=\D\sum_{\scst 0\leq j\leq n}(D_j t^{n-j}\ot T^j\,D_{n-j})D_0\;$  if
$n\geq 0$.

 $\D Aw_4 (t^{\ti 2}Sh_4)^n(t\, Sh_4)^m$\\
\mbox{ } \mbox{ } $\D=\sum_{\scst 0\leq j\leq n}(D_j t^{n-j}\ot T^j\,D_{n-j})D_m
+\sum_{\stackrel{\scst 0\leq j\leq n}{\scst j+n \mbox{ \tiny
is odd}}}(D_{j+1}\ot T^j\,D_{n-j})TD_{m-1}\;\mbox{ if $n,m>0$.}$

 $Aw_4 (t^{\ti 2}Sh_4)^{n_0}(t\, Sh_4)^m (t^{\ti
2}Sh_4)^{n_1}=0\;$ if $n_0,n_1>0$,  $m>1$.

 $Aw_4 (t^{\ti 2}Sh_4)^{n_0}t\, Sh_4 (t^{\ti
2}Sh_4)^{n_1}$\\
\mbox{ } \mbox{ } $=\D \sum_{\stackrel{\scst 0\leq j\leq n_0}{\scst j+n_0 \mbox{ \tiny
is odd}}} (D_{j+1+n_1}\ot T^j\,D_{n_0-j}t^{n_1})TD_{0}\; \mbox{ if $n_0,n_1>0$}
$

 $Aw_4 (t^{\ti 2}Sh_4)^{n_0}
(t\, Sh_4)^{m_1} (t^{\ti 2}Sh_4)^{n_1}\cdots (t\, Sh_4)^{m_k}
 (t^{\ti 2}Sh_4)^{n_{k}} (t\, Sh_4)^m=0\;$\\
\mbox{ } \mbox{ } if $m_i>1$ for some $1\leq i\leq k$ or $n_j$ is even for some $1\leq i<k$.

 $Aw_4 (t^{\ti 2}Sh_4)^{n_0}
t \,Sh_4 (t^{\ti 2}Sh_4)^{n_1}\cdots t \,Sh_4 (t^{\ti 2}Sh_4)^{n_{2k}} (t\, Sh_4)^m$\\
\mbox{ } \mbox{ } $=\left\{\begin{array}{cl}
\D\sum_{\stackrel{\scst 0\leq j\leq n_0}{\scst j+n_0 \mbox{ \tiny
is odd}}} (D_{j+k+\sum n_{2j-1}}t\ot T^j\,D_{-j+k+\sum n_{2j}})D_{m}
&\mbox{ if $n_{2k}$ is even}\\
\mbox{ } \begin{array}{l}
\D\sum_{\stackrel{\scst 0\leq j\leq n_0}{\scst j+n_0 \mbox{ \tiny
is odd}}} (D_{j+k+\sum n_{2j-1}}t\ot T^j\,D_{-j+k+\sum n_{2j}})D_{m} \\
\; \D+\sum_{\stackrel{\scst 0\leq j\leq n_0}{\scst j+n_0 \mbox{ \tiny is odd}}} (D_{j+k+1+\sum n_{2j-1}}\ot
T^j\,D_{-j+k+\sum n_{2j}})TD_{m-1}\end{array}&
\mbox{ otherwise}
\end{array}\right.$\\
\mbox{ } \mbox{ } for $k,n_0,n_1,\dots,n_{2k}>0$, $m\geq 0$ and $n_i$ being odd for $1\leq i\leq 2k-1$.

 $Aw_4 (t^{\ti 2}Sh_4)^{n_0}
t \,Sh_4 (t^{\ti 2}Sh_4)^{n_1}\cdots t\, Sh_4 (t^{\ti 2}Sh_4)^{n_{2k-1}} (t\, Sh_4)^m$\\
\mbox{ } \mbox{ } $
=\left\{\begin{array}{cl}
\D\sum_{\stackrel{\scst 0\leq j\leq n_0}{\scst j+n_0\mbox{ \tiny
is odd}}} (D_{j+k+\sum n_{2j-1}}\ot T^j\,D_{-j+k-1+\sum n_{2j}}t)TD_{m} &\mbox{ if $n_{2k-1}$ is even}\\
\begin{array}{l}
\D\sum_{\stackrel{\scst 0\leq j\leq n_0}{\scst j+n_0\mbox{ \tiny
is odd}}} (D_{j+k+\sum n_{2j-1}}\ot T^j\,D_{-j+k-1+\sum n_{2j}}t)TD_{m} \\
\;\D +\sum_{\stackrel{\scst 0\leq j\leq n_0}{\scst j+n_0
\mbox{ \tiny is odd}}} (D_{j+k+\sum n_{2j-1}}\ot
T^j\,D_{-j+k+\sum n_{2j}})D_{m-1}\end{array}
&\mbox{ otherwise}
\end{array}\right.$\\
\mbox{ } \mbox{ } for $k,n_0,n_1,\dots,n_{2k-1}>0$, $m\geq 0$ and $n_i$ being odd for $1\leq i\leq 2k-2$.

where $D_i=Aw (t \,Sh)^n$.

Let us observe that all the summands  we will obtain when (\ref{pono})
is simplified, has the form:
$\mu c^{\ot 4}(D_a\ot D_b)D_{2i+2-\ell}\Delta(x)$
and
$\mu c^{\ot 4}(D_d\ot D_e)D_{i-k}\Delta(x)$
where $c$ is a $(i+2)$--cocycle, $a+b=i+\ell$, $d+e=2i+2+k$ and
$0\leq\ell, k\leq i$.
Let us recall that if $c$ is an $m$--cocycle and $c'$ is an
$n$--cocycle, then  $c\smile_r c'=\mu (c\ot c')D_r$ is
null if $m< r$ or $n< r$. Taking in mind this last fact, we have
that $\ell+2\geq a$ and $\ell+2\geq b$. Since $a+b=i+\ell$, then
 $i-4\leq\ell\leq i$. On the other hand, $d\leq i+2$
and $e\leq i+2$. Since $d+e=2i+2+k$
 then $0\leq k\leq 2$.
Therefore, the possible non--null summands  of (\ref{pono}) have the form:

$\mu c^{\ot 4}(D_{i+2}\ot D_{i+2})D_{i-2}\Delta(x)$,
 $\mu c^{\ot 4}(D_{i+1}\ot D_{i+2})D_{i-1}\Delta(x)$,\\
$\mu c^{\ot 4}(D_{i+2}\ot D_{i+1})D_{i-1}\Delta(x)$,
 $\mu c^{\ot 4}(D_i\ot D_{i+2})D_i\Delta(x)$,\\
$\mu c^{\ot 4}(D_{i+1}\ot D_{i+1})D_i\Delta(x)$,
$\mu c^{\ot 4}(D_{i+2}\ot D_{i})D_i \Delta(x)$,\\
 $\mu c^{\ot 4}(D_{i+2}\ot D_{i-1})D_{i+1}\Delta(x)$,
$\mu c^{\ot 4}(D_{i+1}\ot D_{i})D_{i+1}\Delta(x)$,\\
$\mu c^{\ot 4}(D_{i}\ot D_{i+1})D_{i+1}\Delta(x)$,
$\mu c^{\ot 4}(D_{i-1}\ot D_{i+2})D_{i+1}\Delta(x)$,\\
 $\mu c^{\ot 4}(D_{i-2}\ot D_{i+2})D_{i+2}\Delta(x)$,
$\mu c^{\ot 4}(D_{i-1}\ot D_{i+1})D_{i+2}\Delta(x)$,\\
$\mu c^{\ot 4}(D_{i}\ot D_{i})D_{i+2}\Delta(x)$,
$\mu c^{\ot 4}(D_{i+1}\ot D_{i-1})D_{i+2}\Delta(x)$,\\
$\mu c^{\ot 4}(D_{i+2}\ot D_{i-2})D_{i+2}\Delta(x)$,
 $\mu c^{\ot 4}(D_{i-2}\ot D_{i+1})D_{i+3}\Delta(x)$,\\
$\mu c^{\ot 4}(D_{i-1}\ot D_{i})D_{i+3}\Delta(x)$,
$\mu c^{\ot 4}(D_{i}\ot D_{i-1})D_{i+3}\Delta(x)$,\\
$\mu c^{\ot 4}(D_{i+1}\ot D_{i-2})D_{i+3}\Delta(x)$,
 $\mu c^{\ot 4}(D_{i-2}\ot D_{i})D_{i+4}\Delta(x)$,\\
$\mu c^{\ot 4}(D_{i-1}\ot D_{i-1})D_{i+4}\Delta(x)$,
$\mu c^{\ot 4}(D_{i}\ot D_{i-2})D_{i+4}\Delta(x)$,\\
 $\mu c^{\ot 4}(D_{i-2}\ot D_{i-1})D_{i+5}\Delta(x)$,
$\mu c^{\ot 4}(D_{i-1}\ot D_{i-2})D_{i+5}\Delta(x)$,\\
 $\mu c^{\ot 4}(D_{i-2}\ot D_{i-2})D_{i+6}\Delta(x)$.

Using the properties mentioned above, the possible non--null summands of (\ref{pono})
for $\ell=i-2$ are:

 $\mu c^{\ot 4} Aw_4 (t^{\ti 2}Sh_4)^{2i+4}(t\,Sh_4)^{i-2}\Delta(x)$
$=\mu c^{\ot 4}(D_{i+2}\ot D_{i+2})D_{i-2}\Delta(x)$,\\
$\mu c^{\ot 4} Aw_4 t^{\ti 2}Sh_4
(t\,Sh_4t^{\ti 2}Sh_4)^{i-3}(t\,Sh_4)^{i+7}\Delta(x)$
$=\mu c^{\ot 4}(D_{i-2}\ot D_{i-2})D_{i+2}\Delta(x)$.

The possible non--null summands of (\ref{pono})
for $\ell=i-1$ are:

$\mu c^{\ot 4} Aw_4(t^{\ti 2}Sh_4)^{2i+3}
(t\,Sh_4)^{i-1}\Delta(x)$\\
\mbox{ } \mbox{ }$=\mu c^{\ot 4}(D_{i+1}\ot D_{i+2}+D_{i+2}\ot D_{i+1})D_{i-1}\Delta(x)$,\\
$\mu c^{\ot 4} Aw_4(t^{\ti 2}Sh_4)^{2i+3-j}t\,Sh_4
(t^{\ti 2}Sh_4)^{j}
(t\,Sh_4)^{i-2}\Delta(x)=0\;$ for $1\leq j\leq 2i+2$,\\
$\mu c^{\ot 4} Aw_4 t^{\ti 2}Sh_4
(t\,Sh_4t^{\ti 2}Sh_4)^{i-2}(t\,Sh_4)^{i+5}\Delta(x)$\\
\mbox{ } \mbox{ }$=\mu c^{\ot 4}((D_{i-2}\ot D_{i-1})D_{i+5}+(D_{i-1}\ot D_{i-1})D_{i+4})
\Delta(x)$,\\
$\mu c^{\ot 4} Aw_4 (t^{\ti 2}Sh_4)^2
(t\,Sh_4t^{\ti 2}Sh_4)^{i-3}(t\,Sh_4)^{i+6}\Delta(x)$\\
\mbox{ } \mbox{ }$=\mu c^{\ot 4}(D_{i-1}\ot D_{i-2})D_{i+5}\Delta(x)$,\\
$\mu c^{\ot 4} Aw_4 t^{\ti 2}Sh_4
(t\,Sh_4t^{\ti 2}Sh_4)^{i-4}t\,Sh_4(t^{\ti 2}Sh_4)^2(t\,Sh_4)^{i+6}\Delta(x)=0$,\\
$\mu c^{\ot 4} Aw_4 (t^{\ti 2}Sh_4)^3
(t\,Sh_4t^{\ti 2}Sh_4)^{i-4}(t\,Sh_4)^{i+7}\Delta(x)=0$,\\
$\mu c^{\ot 4} Aw_4 t^{\ti 2}Sh_4
(t\,Sh_4t^{\ti 2}Sh_4)^{i-5}
t\,Sh_4(t^{\ti 2}Sh_4)^3
(t\,Sh_4)^{i+7}\Delta(x)=0$,\\
$\mu c^{\ot 4} Aw_4 t^{\ti 2}Sh_4
t\,Sh_4(t^{\ti 2}Sh_4)^3
(t\,Sh_4t^{\ti 2}Sh_4)^{i-5}
(t\,Sh_4)^{i+7}\Delta(x)=0$,\\
$\mu c^{\ot 4} Aw_4 t^{\ti 2}Sh_4
(t\,Sh_4t^{\ti 2}Sh_4)^{i-5-j}
t\,Sh_4(t^{\ti 2}Sh_4)^3
(t\,Sh_4t^{\ti 2}Sh_4)^{j}
(t\,Sh_4)^{i+7}\Delta(x)$\\
\mbox{ } \mbox{ }$=0\;$ for $1\leq j\leq i-6$.

Finally, the possible non--null summands of (\ref{pono})
for $\ell=i$ are:

$\mu c^{\ot 4} Aw_4(t^{\ti 2}Sh_4)^{2i+2}(t\,Sh_4)^i\Delta(x)$\\
\mbox{ } \mbox{ }$=\left\{\begin{array}{cl}
\mu c^{\ot 4}((D_i\ot D_{i+2}+D_{i+1}\ot D_{i+1}+D_{i+2}\ot D_i)D_i+
(D_{i+1}\ot D_{i+2})D_{i-1})\Delta(x)& \mbox{ if $i$ is odd}\\
\mu c^{\ot 4}((D_i\ot D_{i+2}+D_{i+1}\ot D_{i+1}+D_{i+2}\ot D_i)D_i+
(D_{i+2}\ot D_{i+1})D_{i-1})\Delta(x)& \mbox{ if $i$ is even},
\end{array}\right.$\\
$\mu c^{\ot 4} Aw_4(t^{\ti 2}Sh_4)^{2i+2-j}t\,Sh (t^{\ti 2}Sh_4)^{j}(t\,Sh_4)^{i-1}\Delta(x)$\\
\mbox{ } \mbox{ }$=\left\{\begin{array}{cl}
\mu c^{\ot 4}(D_{i+1}\ot D_{i+2})D_{i-1}\Delta(x)& \mbox{ if $i$ is odd}\\
\mu c^{\ot 4}(D_{i+2}\ot D_{i+1})D_{i-1}\Delta(x) & \mbox{ if $j,i$ are even}\\
\mu c^{\ot 4}((D_{i+2}\ot D_{i+1})D_{i-1}+D_{i+2}\ot D_{i+2})D_{i-2})\Delta(x) &
\mbox{ if $j$ is odd, $i$ is even}
\end{array}\right.$\\
\mbox{ } \mbox{ }for $1\leq j\leq 2i+1$,\\
$\mu c^{\ot 4} Aw_4(t^{\ti 2}Sh_4)^{2i+2-j-k}t\,Sh (t^{\ti 2}Sh_4)^{j}
t\,Sh (t^{\ti 2}Sh_4)^{k}(t\,Sh_4)^{i-2}\Delta(x)$\\
\mbox{ } \mbox{ }$=\left\{\begin{array}{cl}
\mu c^{\ot 4}(D_{i+2}\ot D_{i+2})D_{i-2}\Delta(x)& \mbox{ if $k,i$ are odd}\\
\mu c^{\ot 4}(D_{i+2}\ot D_{i+2})D_{i-2}\Delta(x) & \mbox{ if $k,i$ are even}\\
0 & \mbox{ otherwise}
\end{array}\right.$\\
\mbox{ } \mbox{ }for $1\leq j,k$, $j+k\leq 2i+1$ and $j$ being odd,\\
$\mu c^{\ot 4} Aw_4 t^{\ti 2}Sh_4(t\,Sh\, t^{\ti 2}Sh_4)^{i-1}(t\,Sh_4)^{i+3}\Delta(x)$\\
\mbox{ } \mbox{ }$=\left\{\begin{array}{cl}
\mu c^{\ot 4}((D_{i-1}\ot D_{i})D_{i+3}+(D_{i}\ot D_{i})D_{i+2})\Delta(x)& \mbox{ if $i$ is odd}\\
\mu c^{\ot 4}((D_{i}\ot D_{i-1})D_{i+3}+(D_{i}\ot D_{i})D_{i+2})\Delta(x) & \mbox{ if $i$ is even},
\end{array}\right.$\\
$\mu c^{\ot 4} Aw_4(t^{\ti 2}Sh_4)^2(t\,Sh\, t^{\ti 2}Sh_4)^{i-2}(t\,Sh_4)^{i+4}\Delta(x)$\\
\mbox{ } \mbox{ }$=\left\{\begin{array}{cl}
\mu c^{\ot 4}((D_{i}\ot D_{i-2})D_{i+4}+(D_{i}\ot D_{i-1})D_{i+3})\Delta(x)& \mbox{ if $i$ is odd}\\
\mu c^{\ot 4}((D_{i-1}\ot D_{i-1})D_{i+4}+(D_{i}\ot D_{i-1})D_{i+3})\Delta(x) & \mbox{ if $i$ is even}
\end{array}\right.$\\
$\mu c^{\ot 4} Aw_4t^{\ti 2}Sh_4 (t\,Sh \,t^{\ti 2}Sh_4)^{i-3}t\,Sh(t^{\ti 2}Sh_4)^2(t\,Sh_4)^{i+4}\Delta(x)$\\
\mbox{ } \mbox{ }$=\left\{\begin{array}{cl}
\mu c^{\ot 4}(D_{i}\ot D_{i-2})D_{i+4}\Delta(x)& \mbox{ if $i$ is odd}\\
\mu c^{\ot 4}(D_{i-2}\ot D_{i})D_{i+4}\Delta(x) & \mbox{ if $i$ is even},
\end{array}\right.$\\
$\mu c^{\ot 4} Aw_4(t^{\ti 2}Sh_4)^3 (t\,Sh\, t^{\ti 2}Sh_4)^{i-3}(t\,Sh_4)^{i+5}\Delta(x)$\\
\mbox{ } \mbox{ }$=\left\{\begin{array}{cl}
\mu c^{\ot 4}((D_{i-1}\ot D_{i-2})D_{i+5}+(D_{i-2}\ot D_{i}+D_{i}\ot D_{i-2})D_{i+4})\Delta(x)
& \mbox{ if $i$ is odd}\\
\mu c^{\ot 4}((D_{i-2}\ot D_{i-1})D_{i+5}+(D_{i-2}\ot D_{i}+D_{i}\ot D_{i-2})D_{i+4})\Delta(x)
& \mbox{ if $i$ is even},
\end{array}\right.$\\
$\mu c^{\ot 4} Aw_4(t^{\ti 2}Sh_4)^2 (t\,Sh\, t^{\ti 2}Sh_4)^{i-4}t\,Sh(t^{\ti 2}Sh_4)^2
(t\,Sh_4)^{i+5}\Delta(x)$\\
\mbox{ } \mbox{ }$=\left\{\begin{array}{cl}
\mu c^{\ot 4}(D_{i-2}\ot D_{i-1})D_{i+5}\Delta(x)
& \mbox{ if $i$ is odd}\\
0
& \mbox{ if $i$ is even},
\end{array}\right.$\\
$\mu c^{\ot 4} Aw_4t^{\ti 2}Sh_4 (t\,Sh\, t^{\ti 2}Sh_4)^{i-4}t\,Sh(t^{\ti 2}Sh_4)^3
(t\,Sh_4)^{i+5}\Delta(x)$\\
\mbox{ } \mbox{ }$=\left\{\begin{array}{cl}
\mu c^{\ot 4}(D_{i-2}\ot D_{i})D_{i+4}\Delta(x)
& \mbox{ if $i$ is odd}\\
\mu c^{\ot 4}(D_{i}\ot D_{i-2})D_{i+4}\Delta(x)
& \mbox{ if $i$ is even},
\end{array}\right.$\\
$\mu c^{\ot 4} Aw_4t^{\ti 2}Sh_4 t\,Sh (t^{\ti 2}Sh_4)^{3}
(t\,Sh\, t^{\ti 2}Sh_4)^{i-4}(t\,Sh_4)^{i+5}\Delta(x)$\\
\mbox{ } \mbox{ }$=\left\{\begin{array}{cl}
\mu c^{\ot 4}((D_{i-1}\ot D_{i-2})D_{i+5}+(D_{i}\ot D_{i-2})D_{i+4})\Delta(x)
& \mbox{ if $i$ is odd}\\
\mu c^{\ot 4}(D_{i}\ot D_{i-2})D_{i+4}\Delta(x)
& \mbox{ if $i$ is even},
\end{array}\right.$\\
$\mu c^{\ot 4} Aw_4t^{\ti 2}Sh_4
(t\,Sh\, t^{\ti 2}Sh_4)^{i-4-j}
t\,Sh (t^{\ti 2}Sh_4)^{3}
(t\,Sh\, t^{\ti 2}Sh_4)^{j}(t\,Sh_4)^{i+5}\Delta(x)$\\
\mbox{ } \mbox{ }$=\left\{\begin{array}{cl}
\mu c^{\ot 4}((D_{i-1}\ot D_{i-2})D_{i+5}+(D_{i}\ot D_{i-2})D_{i+4})\Delta(x)
& \mbox{ if $j,i$ are odd}\\
\mu c^{\ot 4}(D_{i-2}\ot D_{i})D_{i+4}\Delta(x)
& \mbox{ if $j$ is even, $i$ is odd}\\
\mu c^{\ot 4}((D_{i-2}\ot D_{i-1})D_{i+5}+(D_{i-2}\ot D_{i})D_{i+4})\Delta(x)
& \mbox{ if $j$ is odd, $i$ is even}\\
\mu c^{\ot 4}(D_{i}\ot D_{i-2})D_{i+4}\Delta(x)
& \mbox{ if $j,i$ is even}
\end{array}\right.$\\
\mbox{ } \mbox{ }for $1\leq j\leq i-5$,\\
$\mu c^{\ot 4} Aw_4(t^{\ti 2}Sh_4)^4
(t\,Sh\, t^{\ti 2}Sh_4)^{i-4}(t\,Sh_4)^{i+6}\Delta(x)$\\
\mbox{ } \mbox{ }$=\left\{\begin{array}{cl}
\mu c^{\ot 4}((D_{i-2}\ot D_{i-2})D_{i+6}+(D_{i-2}\ot D_{i-1})D_{i+5})\Delta(x)
& \mbox{ if $i$ is odd}\\
\mu c^{\ot 4}(D_{i-2}\ot D_{i-1})D_{i+5}\Delta(x)
& \mbox{ if $i$ is even},
\end{array}\right.$\\
$\mu c^{\ot 4} Aw_4(t^{\ti 2}Sh_4)^3
(t\,Sh\, t^{\ti 2}Sh_4)^{i-5}t\,Sh (t^{\ti 2}Sh_4)^{2}
(t\,Sh_4)^{i+6}\Delta(x)$\\
\mbox{ } \mbox{ }$=\mu c^{\ot 4}(D_{i-2}\ot D_{i-2})D_{i+6}\Delta(x)$,\\
$\mu c^{\ot 4} Aw_4(t^{\ti 2}Sh_4)^2
(t\,Sh\, t^{\ti 2}Sh_4)^{i-5}t\,Sh (t^{\ti 2}Sh_4)^{3}
(t\,Sh_4)^{i+6}\Delta(x)$\\
\mbox{ } \mbox{ }$=\left\{\begin{array}{cl}
0
& \mbox{ if $i$ is odd}\\
\mu c^{\ot 4}(D_{i-2}\ot D_{i-1})D_{i+5}\Delta(x)
& \mbox{ if $i$ is even},
\end{array}\right.$\\
$\mu c^{\ot 4} Aw_4(t^{\ti 2}Sh_4)^2 t\,Sh (t^{\ti 2}Sh_4)^{3}
(t\,Sh\, t^{\ti 2}Sh_4)^{i-5}
(t\,Sh_4)^{i+6}\Delta(x)=0$,\\
$\mu c^{\ot 4} Aw_4(t^{\ti 2}Sh_4)^2
(t\,Sh\, t^{\ti 2}Sh_4)^{i-5-j}
t\,Sh (t^{\ti 2}Sh_4)^{3}
(t\,Sh\, t^{\ti 2}Sh_4)^{j}
(t\,Sh_4)^{i+6}\Delta(x)$\\
\mbox{ } \mbox{ }$=\left\{\begin{array}{cl}
\mu c^{\ot 4}((D_{i-2}\ot D_{i-2})D_{i+6}+(D_{i-2}\ot D_{i-1})D_{i+5})\Delta(x)
& \mbox{ if $j,i$ are odd}\\
\mu c^{\ot 4}(D_{i-2}\ot D_{i-1})D_{i+5}\Delta(x)
& \mbox{ if $j,i$ are even}\\
0& \mbox{ otherwise}
\end{array}\right.$\\
\mbox{ } \mbox{ }for $1\leq j\leq i-6$,\\
$\mu c^{\ot 4} Aw_4t^{\ti 2}Sh_4
(t\,Sh\, t^{\ti 2}Sh_4)^{i-5}
t\,Sh (t^{\ti 2}Sh_4)^{4}
(t\,Sh_4)^{i+6}\Delta(x)=0$,\\
$\mu c^{\ot 4} Aw_4t^{\ti 2}Sh_4
t\,Sh (t^{\ti 2}Sh_4)^{3}
(t\,Sh\, t^{\ti 2}Sh_4)^{i-6}
t\,Sh (t^{\ti 2}Sh_4)^{2}
(t\,Sh_4)^{i+6}\Delta(x)$\\
\mbox{ } \mbox{ }$=\left\{\begin{array}{cl}
\mu c^{\ot 4}(D_{i-2}\ot D_{i-2})D_{i+6}\Delta(x)
& \mbox{ if $i$ is odd}\\
0
& \mbox{ if $i$ is even},
\end{array}\right.$\\
$\mu c^{\ot 4} Aw_4t^{\ti 2}Sh_4
(t\,Sh\, t^{\ti 2}Sh_4)^{i-6-j}
t\,Sh (t^{\ti 2}Sh_4)^{3}
(t\,Sh\, t^{\ti 2}Sh_4)^{j}
t\,Sh (t^{\ti 2}Sh_4)^{2}
(t\,Sh_4)^{i+6}\Delta(x)$\\
\mbox{ } \mbox{ }$=\left\{\begin{array}{cl}
\mu c^{\ot 4}(D_{i-2}\ot D_{i-2})D_{i+6}\Delta(x)
& \mbox{ if $j$ is even}\\
0&\mbox{ otherwise}
\end{array}\right.$\\
\mbox{ } \mbox{ }for $1\leq j\leq i-7$,\\
$\mu c^{\ot 4} Aw_4t^{\ti 2}Sh_4
(t\,Sh\, t^{\ti 2}Sh_4)^{i-6}
t\,Sh (t^{\ti 2}Sh_4)^{3}
t\,Sh (t^{\ti 2}Sh_4)^{2}
(t\,Sh_4)^{i+6}\Delta(x)$\\
\mbox{ } \mbox{ }$=\mu c^{\ot 4}(D_{i-2}\ot D_{i-2})D_{i+6}\Delta(x)$,\\
$\mu c^{\ot 4} Aw_4(t^{\ti 2}Sh_4)^5
(t\,Sh\, t^{\ti 2}Sh_4)^{i-5}
(t\,Sh_4)^{i+7}\Delta(x)$\\
\mbox{ } \mbox{ }$=\mu c^{\ot 4}(D_{i-2}\ot D_{i-2})D_{i+6}\Delta(x)$,\\
$\mu c^{\ot 4} Aw_4t^{\ti 2}Sh_4
(t\,Sh\, t^{\ti 2}Sh_4)^{i-6}
t\,Sh (t^{\ti 2}Sh_4)^{5}
(t\,Sh_4)^{i+7}\Delta(x)=0$,\\
$\mu c^{\ot 4} Aw_4t^{\ti 2}Sh_4
(t\,Sh\, t^{\ti 2}Sh_4)^{i-6-j}
t\,Sh (t^{\ti 2}Sh_4)^{5}
(t\,Sh\, t^{\ti 2}Sh_4)^{j}
(t\,Sh_4)^{i+7}\Delta(x)$\\
\mbox{ } \mbox{ }$=0\;$ for $1\leq j\leq i-7$,\\
$\mu c^{\ot 4} Aw_4(t^{\ti 2}Sh_4)^3
(t\,Sh\, t^{\ti 2}Sh_4)^{i-6}
t\,Sh (t^{\ti 2}Sh_4)^{3}
(t\,Sh_4)^{i+7}\Delta(x)$\\
\mbox{ } \mbox{ }$=\mu c^{\ot 4}(D_{i-2}\ot D_{i-2})D_{i+6}\Delta(x)$,\\
$\mu c^{\ot 4} Aw_4(t^{\ti 2}Sh_4)^3
t\,Sh (t^{\ti 2}Sh_4)^{3}
(t\,Sh\, t^{\ti 2}Sh_4)^{i-6}
(t\,Sh_4)^{i+7}\Delta(x)$\\
\mbox{ } \mbox{ }$=\mu c^{\ot 4}(D_{i-2}\ot D_{i-2})D_{i+6}\Delta(x)$,\\
$\mu c^{\ot 4} Aw_4(t^{\ti 2}Sh_4)^3
(t\,Sh\, t^{\ti 2}Sh_4)^{i-6-j}
t\,Sh (t^{\ti 2}Sh_4)^{3}
(t\,Sh\, t^{\ti 2}Sh_4)^{j}
(t\,Sh_4)^{i+7}\Delta(x)$\\
\mbox{ } \mbox{ }$=\mu c^{\ot 4}(D_{i-2}\ot D_{i-2})D_{i+6}\Delta(x)\;$
for $1\leq j\leq i-7$,\\
$\mu c^{\ot 4} Aw_4t^{\ti 2}Sh_4
t\,Sh (t^{\ti 2}Sh_4)^{3}
(t\,Sh\, t^{\ti 2}Sh_4)^{i-7}
t\,Sh (t^{\ti 2}Sh_4)^{3}
(t\,Sh_4)^{i+7}\Delta(x)$\\
\mbox{ } \mbox{ }$=\left\{\begin{array}{cl}
\mu c^{\ot 4}(D_{i-2}\ot D_{i-2})D_{i+6}\Delta(x)
& \mbox{ if $i$ is odd}\\
0&\mbox{ if $i$ is even},
\end{array}\right.$\\
$\mu c^{\ot 4} Aw_4t^{\ti 2}Sh_4
(t\,Sh\, t^{\ti 2}Sh_4)^{i-7}
t\,Sh (t^{\ti 2}Sh_4)^{3}
t\,Sh (t^{\ti 2}Sh_4)^{3}
(t\,Sh_4)^{i+7}\Delta(x)$\\
\mbox{ } \mbox{ }$=\mu c^{\ot 4}(D_{i-2}\ot D_{i-2})D_{i+6}\Delta(x)$,\\
$\mu c^{\ot 4} Aw_4t^{\ti 2}Sh_4
(t\,Sh\, t^{\ti 2}Sh_4)^{i-7-j}
t\,Sh (t^{\ti 2}Sh_4)^{3}
(t\,Sh\, t^{\ti 2}Sh_4)^{j}
t\,Sh (t^{\ti 2}Sh_4)^{3}
(t\,Sh_4)^{i+7}\Delta(x)$,\\
\mbox{ } \mbox{ }$=\left\{\begin{array}{cl}
\mu c^{\ot 4}(D_{i-2}\ot D_{i-2})D_{i+6}\Delta(x)
& \mbox{ if $j$ is even}\\
0&\mbox{ otherwise}
\end{array}\right.$\\
\mbox{ } \mbox{ }for $1\leq j\leq i-8$.

Now, using property (\ref{pie}),
 the simplified expression for (\ref{pono}) mod $2$ is:

$\mu c^{\ot 4}((D_0\ot D_{2}+D_{2}\ot D_0+D_{1}\ot D_{1})D_0
+(D_0\ot D_0)D_{2})\Delta(x)$\\
\mbox{ } \mbox{ }$=\delta(c\smile_0 c\smile_1 c\smile_2 c)
+c\smile_1 c\smile_0 c\smile_1 c+c\smile_0 c\smile_2 c\smile_0 c\;$
 for $i=0$.

$\mu c^{\ot 4}((D_1\ot D_{3}+D_{3}\ot D_1+D_{2}\ot D_{2})D_1+(D_1\ot D_1)D_{3})\Delta(x)$\\
\mbox{ } \mbox{ }$=\delta(c\smile_1 c\smile_2 c\smile_3 c)+c\smile_2 c\smile_1 c\smile_2 c
+c\smile_1 c\smile_3 c\smile_1 c\;$ for $i=1$.

$\mu c^{\ot 4}((D_2\ot D_{4}+D_{4}\ot D_2+D_{3}\ot D_{3})D_2+(D_2\ot D_2)D_{4}$\\
\mbox{ } \mbox{ }$\quad+(D_{0}\ot D_2+D_2\ot D_{0})D_{6})\Delta(x)$\\
\mbox{ } \mbox{ }$=\delta(c\smile_2 c\smile_3 c\smile_4 c+c\smile_0 c\smile_7 c\smile_2 c)
+c\smile_3 c\smile_2 c\smile_3 c$\\
\mbox{ } \mbox{ }$\quad+c\smile_2 c\smile_4 c\smile_2 c\;$ for $i=2$.

 $\mu c^{\ot 4}((D_3\ot D_{5}+D_{5}\ot D_3+D_{4}\ot D_{4})D_3+(D_3\ot D_3)D_{5}$\\
\mbox{ } \mbox{ }$\quad+(D_{2}\ot D_{3}+D_3\ot D_{2})D_{6}+
(D_{3}\ot D_{1}+D_{1}\ot D_{3})D_{7}
)\Delta(x)$\\
\mbox{ } \mbox{ }$=\delta(c\smile_3 c\smile_4 c\smile_5 c+ c\smile_2 c\smile_7 c\smile_3 c
+c\smile_3 c\smile_8 c\smile_1 c)$\\
\mbox{ } \mbox{ }$\quad+c\smile_3 c\smile_5 c\smile_3 c
+c\smile_4 c\smile_3 c\smile_4 c\;$ for $i=3$.

$\mu c^{\ot 4}((D_4\ot D_{6}+D_{6}\ot D_4+D_{5}\ot D_{5})D_4+(D_4\ot D_4)D_{6})\Delta(x)$\\
\mbox{ } \mbox{ }$=\delta(c\smile_4 c\smile_5 c\smile_6 c)
+c\smile_5 c\smile_4 c\smile_5 c+c\smile_4 c\smile_6 c\smile_4 c\;$ for $i=4$.

$\mu c^{\ot 4}((D_5\ot D_{7}+D_{7}\ot D_5+D_{6}\ot D_{6})D_5+(D_5\ot D_5)D_{7}$\\
\mbox{ } \mbox{ }$\quad+(D_{4}\ot D_{5}+D_5\ot D_{4})D_{8}+
(D_{4}\ot D_{3}+D_{3}\ot D_{4})D_{10})\Delta(x)$\\
\mbox{ } \mbox{ }$=\delta(c\smile_5 c\smile_6 c\smile_7 c+ c\smile_4 c\smile_9 c\smile_5 c
+c\smile_4 c\smile_{11} c\smile_3 c)$\\
\mbox{ } \mbox{ }$\quad +c\smile_5 c\smile_7 c\smile_5 c
+c\smile_6 c\smile_5 c\smile_6 c\;$  for $i=5$.

$\mu c^{\ot 4}((D_i\ot D_{i+2}+D_{i+2}\ot D_i
+D_{i+1}\ot D_{i+1})D_i+(D_i\ot D_i)D_{i+2})\Delta(x)$\\
\mbox{ } \mbox{ }$=\delta(c\smile_i c\smile_{i+1} c\smile_{i+2} c)
+c\smile_{i+1} c\smile_i c\smile_{i+1} c
+c\smile_i c\smile_{i+2} c\smile_i c$\\
\mbox{ } \mbox{ } for $i\geq 6$,  $i$ being even, $\frac{i-6}{2}$ being odd.

 $\mu c^{\ot 4}((D_i\ot D_{i+2}+D_{i+2}\ot D_i+D_{i+1}\ot D_{i+1})D_i
+(D_i\ot D_i)D_{i+2}$\\
\mbox{ } \mbox{ }$\quad+(D_i\ot D_{i-2}+D_{i-2}\ot D_i)D_{i+4}$\\
\mbox{ } \mbox{ }$=\delta(c\smile_i c\smile_{i+1} c\smile_{i+2} c+
c\smile_i c\smile_{i+5} c\smile_{i-2} c)$\\
\mbox{ } \mbox{ }$\quad+c\smile_{i+1} c\smile_i c\smile_{i+1} c
+c\smile_i c\smile_{i+2} c\smile_i c$\\
\mbox{ } \mbox{ } for $i\geq 6$ and $i,\frac{i-6}{2}$ being even.

$\mu c^{\ot 4}((D_i\ot D_{i+2}+D_{i+2}\ot D_i+D_{i+1}\ot D_{i+1})D_i
+(D_i\ot D_i)D_{i+2}$\\
\mbox{ } \mbox{ }$\quad+(D_{i-1}\ot D_{i}+D_i\ot D_{i-1})D_{i+3}
+(D_{i}\ot D_{i-2}+D_{i-2}\ot D_{i})D_{i+4}$\\
\mbox{ } \mbox{ }$\quad+(D_{i-1}\ot D_{i-2}+D_{i-2}\ot D_{i-1})D_{i+5}
)\Delta(x)$\\
\mbox{ } \mbox{ }$=\delta(c\smile_i c\smile_{i+1} c\smile_{i+2} c+ c\smile_i c\smile_{i+4} c\smile_{i-1} c
+c\smile_i c\smile_{i+5} c\smile_{i-2} c$\\
\mbox{ } \mbox{ }$\quad+c\smile_{i-1} c\smile_{i+6} c\smile_{i-2} c)
+c\smile_{i+1} c\smile_i c\smile_{i+1} c
+c\smile_i c\smile_{i+2} c\smile_i c$\\
 \mbox{ } \mbox{ } for  $i\geq 7$ and $i,\frac{i-7}{2}$ being odd.

$\mu c^{\ot 4}(
(D_i\ot D_{i+2}+D_{i+2}\ot D_i+D_{i+1}\ot D_{i+1})D_i
+(D_i\ot D_i)D_{i+2} $\\
\mbox{ } \mbox{ }$\quad+(D_{i-1}\ot D_{i}+D_i\ot D_{i-1})D_{i+3}
+(D_{i}\ot D_{i-2}+D_{i-2}\ot D_{i})D_{i+4}$\\
\mbox{ } \mbox{ }$\quad+(D_{i-1}\ot D_{i-2}+D_{i-2}\ot D_{i-1})D_{i+5}
)\Delta(x)$\\
\mbox{ } \mbox{ }$=\delta(c\smile_i c\smile_{i+1} c\smile_{i+2} c+ c\smile_i c\smile_{i+4} c\smile_{i-1} c
+c\smile_i c\smile_{i+5} c\smile_{i-2} c$\\
\mbox{ } \mbox{ }$\quad+c\smile_{i-1} c\smile_{i+6} c\smile_{i-2} c)
+c\smile_{i+1} c\smile_i c\smile_{i+1} c+c\smile_i c\smile_{i+2} c\smile_i c$\\
\mbox{ } \mbox{ } for  $i\geq 7$, $i$ being odd, $\frac{i-7}{2}$ being even.

This completes the proof.
\end{document}